\documentclass{amsart}

\usepackage{amsfonts,amssymb,amsmath,amscd,amsthm,color}
\usepackage[all]{xy}

\numberwithin{equation}{section}

\newtheorem{theorem}{Theorem}[section]
\newtheorem{lemma}[theorem]{Lemma}
\newtheorem{proposition}[theorem]{Proposition}
\newtheorem{corollary}[theorem]{Corollary}
\newtheorem{conjecture}[theorem]{Conjecture}

\theoremstyle{definition}
\newtheorem{definition}[theorem]{Definition}

\theoremstyle{remark}
\newtheorem{claim}[theorem]{Claim}

\newtheorem{notation}{Notation}
\newtheorem{remark}{Remark}[section]
\newtheorem{example}{Example}[section]

\newcommand{\Hom}{{\rm Hom}}

\newcommand{\Spec}{{\rm Spec}}

\newcommand{\CCC}{{\mathbb C}}
\newcommand{\FF}{{\mathbb F}}
\newcommand{\EE}{{\mathbb E}}
\newcommand{\RR}{{\mathbb R}}
\newcommand{\ZZ}{{\mathbb Z}}
\newcommand{\PP}{{\mathbb P}}
\newcommand{\NN}{{\mathbb N}}
\newcommand{\GG}{{\mathbb G}}

\newcommand{\QQ}{{\mathbb Q}}
\newcommand{\KK}{{\mathbb K}}
\newcommand{\LL}{{\mathbb L}}

\newcommand{\AAA}{{\mathbb A}}

\newcommand{\kk}{{\Bbbk}}
\newcommand{\CL}{{\mathcal L}}
\newcommand{\CA}{{\mathcal A}}
\newcommand{\CO}{{\mathcal O}}

\newcommand{\CN}{{\mathcal N}}

\newcommand{\CF}{{\mathcal F}}

\newcommand{\CE}{{\mathcal E}}

\begin{document}
\author{Ilya Tyomkin}
\title{On Zariski's theorem in positive characteristic}
\address{Department of Mathematics, Ben-Gurion University of the Negev, P.O.Box 653, Be'er Sheva, 84105, Israel}
\email{tyomkin@cs.bgu.ac.il}
\thanks{The research leading to these results has received funding from the European Union Seventh Framework Programme (FP7/2007-2013) under grant agreement 248826.}
\begin{abstract}
In the current paper we show that the dimension of a family $V$ of irreducible reduced curves in a given ample linear system on a toric surface $S$ over an algebraically closed field is bounded from above by $-K_S.C+p_g(C)-1$, where $C$ denotes a general curve in the family. This result generalizes a famous theorem of Zariski to the case of positive characteristic. We also explore new phenomena that occur in positive characteristic: We show that the equality $\dim(V)=-K_S.C+p_g(C)-1$ does not imply the nodality of $C$ even if $C$ belongs to the smooth locus of $S$, and construct reducible Severi varieties on weighted projective planes in positive characteristic, parameterizing irreducible reduced curves of given geometric genus in a given very ample linear system.
\end{abstract}
\maketitle

\section{Introduction}
In 1982, Zariski proved the following remarkable theorem over an algebraically closed field of characteristic zero:
\begin{theorem}[{\cite[Theorem 2]{Zar82}}]
Let $S=\PP^2$ be the projective plane, $\CL\in Pic(S)$ be a line bundle,  $V\subset |\CL|$ be an irreducible subvariety, whose general closed point corresponds to a reduced curve $C$. Then
\begin{enumerate}
\item $\dim(V)\le -K_S.C+p_g(C)-1$; and
\item If the equality holds, then $C$ is nodal.
\end{enumerate}
\end{theorem}
Under certain numerical conditions on $\CL$ and $K_S$, Zariski's theorem was generalized to the case of families of curves on rational surfaces satisfying tangency conditions by Harris \cite[Proposition 2.1]{Har86}, Caporaso-Harris \cite[Propositions~2.1, 2.2]{CH98g}, \cite[Proposition~2.1]{CH98r}, Vakil \cite[Theorem 3.1]{Vak00}, and others. Zariski's theorem and its generalizations played an important role in Harris's proof of the irreducibility of Severi varieties \cite{Har86}, and in a series of enumerative results \cite{CH98g,CH98r,Vak00}.

Two different approaches were developed to prove Zariski's theorem and its  generalizations, and both of them used the assumption on the characteristic. Thus, the question whether Zariski's theorem and its  generalizations hold true in positive characteristic remained open.

In the original Zariski's approach, the assumption was used in order to find a local parameter that would parameterize a general one-dimensional subfamily in $V$ {\em and} the nodes of the corresponding curves \cite[p.216]{Zar82}. Note that in positive characteristic such local parameter does not exist in general.

The second approach was developed by Caporaso and Harris, and is based on a result of Arbarello and Cornalba \cite{AC81}. The main idea was to embed the tangent space $T_xV$, where $x\in V$ is a general closed point, into the space of first-order deformations  $Def^1(C,f)$, where $C$ is the normalization of the curve corresponding to $x$ and $f\colon C\to S$ is the natural map; and to describe the image of $T_xV$ in $Def^1(C,f)$. Thus, the assumption on the characteristic was essential in order to pass from the embedded deformations of curves to deformations of maps from their normalizations to the surface \cite[p.356]{CH98g}, \cite[p.159]{CH98r}. The latter is possible in characteristic zero, since if $C\to T$ is a one-parameter family of curves and $\widetilde{C}\to C$ is its normalization then over an open subset $T'\subset T$, the fibers of $\widetilde{C}\to T$ are smooth, hence coincide with the normalizations of the fibers of $C\to T$, which is no longer true in positive characteristic. In positive characteristic, one may need, first, to proceed with a purely inseparable base change, which destroys the argument based on the first-order computations.

The generalizations of Zariski's theorem can be summarized as follows: Let $S$ be a rational surface, $E\subset S$ be a reduced curve, and $\CL\in Pic(S)$ be a line bundle. For an irreducible component $E'\subset E$, let $\alpha^{E'}=(\alpha^{E'}_1, \alpha^{E'}_2,\dotsc)$ and $\beta^{E'}=(\beta^{E'}_1, \beta^{E'}_2,\dotsc)$ be sequences of non-negative integers such that $\sum_i i\alpha^{E'}_i+\sum_i i\beta^{E'}_i=\CL.E'$, and let $\Omega^{E'}=\{p_{i,j}^{E'}\}_{i,\, 1\le j\le \alpha_i}\subset E'$ be a family of general points. For a sequence of integers $\gamma=(\gamma_1,\gamma_2,\dotsc)$, set $|\gamma|:=\sum_i\gamma_i$ and $I\gamma:=\sum_ii\gamma_i$. Set $|\alpha|:=\sum_{E'}|\alpha^{E'}|$, $|\beta|:=\sum_{E'}|\beta^{E'}|$, and $\Omega:=\cup \Omega^{E'}$; and let $G\subset S$ be any curve disjoint from $\Omega$.
\begin{theorem}\label{thm:genZar}
Let $V\subset |\CL|$ be a positive-dimensional irreducible subvariety, whose general closed point corresponds to a reduced irreducible curve $C$. Assume that $-C.(K_S+E)+|\beta|>1$, $C$ belongs to the smooth locus of $S$, and for any irreducible component $E'\subset E$, the normalization $\widetilde{C}$ contains points $\{q_{i,j}^{E'}\}_{i,\, 1\le j\le \alpha_i}$ and $\{r_{i,j}^{E'}\}_{i,\, 1\le j\le \beta_i}$ such that $q_{i,j}^{E'}$ are mapped to $p_{i,j}^{E'}$, and $f^*(E')=\sum i(q_{i,j}^{E'}+r_{i,j}^{E'})$, where $f\colon\widetilde{C}\to C\hookrightarrow S$ denotes the natural map. Then
\begin{enumerate}
\item $\dim(V)\le -C.(K_S+E)+|\beta|+p_g(C)-1$; and
\item If the equality holds, then $df$ is nowhere zero, $\Omega$ is the set of base points of $V$, $C$ is smooth along its intersection with $E$, and $C$ intersects $G$ transversally. If, in addition, $C$ is singular and $-C.(K_S+E)+|\beta|>3$ then $C$ is nodal.
\end{enumerate}
\end{theorem}

In the current paper we consider the case of pairs $(S,E)$, where $S$ is a toric surfaces and $E=\partial S$ is the complement of the maximal orbit. We give a characteristic-independent proof of the first part of Theorem~\ref{thm:genZar} for such pairs $(S,E)$. We also construct counterexamples to the second part of the theorem in positive characteristic. Our examples include families of plane curves satisfying tangency conditions, and families of curves on weighted projective planes with no tangency conditions imposed. We note here, that so far we have found no counterexamples to the second part of the original Zariski's theorem. Finally, we use these examples to construct {\em reducible} Severi varieties on weighted projective planes in positive characteristic.

Our approach to Theorem~\ref{thm:genZar} (1) is based on the canonical tropicalization procedure in positive characteristic developed in \cite{T10}: We replace the ground field by a field equipped with a non-Archimedean valuation. To a given curve $C\subset S$, we associate a parameterized tropical curve $(\Gamma, h_\Gamma)$ of genus at most $p_g(C)$ in a canonical way. Then we show that the dimension of $V$ is bounded by the dimension of the space of tropical curves satisfying certain conditions. Note, that the latter dimension does not depend on the characteristic! Finally, we show that it is bounded by $-C.(K_S+E)+|\beta|+p_g(C)-1$.

The idea beyond our counterexamples to the second part of Theorem~\ref{thm:genZar} is the following: If the toric surface $S$ is the quotient of another surface $\widetilde{S}$ by an action of $\mu_p$, where $p>2$ is the characteristic, then $\widetilde{S}\to S$ is bijective, and the images of smooth curves $C\subset \widetilde{S}$ may have only unibranch singularities of type $A_{p-1}$. This observation leads to Theorems~\ref{thm:char>2} and \ref{thm:chararb} that show existence of maximal-dimensional families of curves with singularities of type $A_{p-1}$, which contradicts the statement of Theorem~\ref{thm:genZar} (2) in positive characteristic.

In order to construct reducible Severi varieties on toric surfaces in positive characteristic, we use the examples above and some deformation theory to exhibit non-empty components of the Severi varieties, whose general closed points correspond to curves having different types of singularities.

Finally, we would like to recall our conjecture \cite[Conjecture 1.2]{T07}:
\begin{conjecture}
If $S$ is a toric surface over an algebraically closed field of characteristic zero, $\CL\in Pic(S)$ is an effective class, and $g\ge 0$ is a non-negative integer then the Severi variety $V^{irr}(S,\CL,g)$ parameterizing irreducible nodal curves of genus $g$ in the linear system $|\CL|$ that do not contain the zero-dimensional orbits of $S$ is either empty or irreducible.
\end{conjecture}
The conjecture is known to be true in the plane case due to the famous result of Harris \cite{Har86}, and in the case of Hirzebruch surfaces \cite{T07}. In \cite{T07}, we also prove the conjecture for rational curves on any toric surface in arbitrary characteristic. As our examples show, the analog of the conjecture fails in positive characteristic. However, we believe that it does hold true in characteristic zero.

\subsection{Plan of the paper}
In Section~\ref{sec:trop}, we recall the definitions and basic properties of parameterized tropical curves, and of the canonical tropicalization procedure. Section~\ref{sec:mainthmproof} is devoted to the proof of Theorem~\ref{thm:genZar} (1) for toric pairs $(S,E)$ in arbitrary characteristic. In Section~\ref{sec:exandcex}, we construct counterexamples to Theorem~\ref{thm:genZar} (2) (Theorems~\ref{thm:char>2} and \ref{thm:chararb}), and examples of reducible Severi varieties on weighted projective planes (Theorem~\ref{thm:redSev}) in positive characteristic. The necessary deformation theory for the proof of Theorem~\ref{thm:redSev} is discussed in Subsection~\ref{subsec:def}

\subsection{Conventions and notation}\label{subsec:convnot}
\subsubsection{Non-Archimedean base field}
Throughout this paper, $\kk$ denotes an algebraically closed field, $R$ denotes a complete discrete valuation ring with residue field $\kk$ and field of fractions $\FF$, $\overline{\FF}$ denotes the separable closure of $\FF$, and $\nu$ denotes the valuation on $\overline{\FF}$ normalized such that $\nu(\FF^*)=\ZZ$. For an intermediate extension $\FF\subseteq\LL\subseteq \overline{\FF}$, $R_\LL$ denotes the ring of integers in $\LL$. Note that if $[\LL:\FF]<\infty$ then $R_\LL$ is a complete discrete valuation ring since $R$ is so. For two finite intermediate extensions $\FF\subseteq\KK\subseteq\LL\subseteq \overline{\FF}$, the relative ramification index $[\nu(\LL^*):\nu(\KK^*)]$ is denoted by $e_{\LL/\KK}$, and if $\KK=\FF$ then it is denoted simply by $e_{\LL}$. For a finite intermediate extension $\FF\subseteq\LL\subseteq \overline{\FF}$, $t_\LL$ denotes a uniformizer in $R_\LL$.

\subsubsection{Algebraic varieties}
For an algebraic variety $X$ defined over a ring $A$, and an $A$-algebra $B$, we denote by $X(B)$ the set of $B$-points of $X$. If $\CL$ is a line bundle on $X$ then $\CL^d:=\CL^{\otimes d}$ denotes the $d$-th tensor power of $\CL$. If $D\subset X$ is a reduced divisor whose generic points belong to the regular locus of $X$ then $\Omega_X({\rm log}(D))$ denotes the corresponding log-differential forms, i.e., forms having at most simple poles at the generic points of $D$.

\subsubsection{Latices and toric varieties}
In this paper, $M$ denotes a lattice of finite rank, $N:=\Hom_\ZZ(M,\ZZ)$ denotes the dual lattice, and $\Sigma$ denotes a fan in $N_\RR:=N\otimes_\ZZ \RR$. We set $T_N:=\Spec\, \ZZ[M]$ and $N_\QQ:=N\otimes_\ZZ \QQ$. The monomials in $\ZZ[M]$ are denoted by $x^m$. For $\sigma,\tau\in\Sigma$, set $X_\sigma:=\Spec\,\ZZ[\check{\sigma}\cap M]$, $X_{\sigma\tau}:=X_\sigma\cap X_\tau=X_{\sigma\cap\tau}$, and $X_\Sigma:=\cup_{\sigma\in\Sigma} X_\sigma$; and denote by $\partial X_\Sigma$ the complement of $T_N$ in $X_\Sigma$. Recall that $m\mapsto \frac{dx^m}{x^m}$ gives a canonical isomorphism from $M\otimes_\ZZ\CO_X$ to the sheaf of log-differential forms $\Omega_{X_\Sigma}({\rm log} (\partial X_\Sigma))$. Rays, i.e., one-dimensional cones, in $\Sigma$ are denoted by $\rho$. For any ray $\rho\in\Sigma$, the closure of the corresponding codimension-one orbit is denoted by $E_\rho$. Then $\partial X_\Sigma:=\cup_\rho E_\rho$. For a lattice polytope $\Delta\subset M_\RR:=M\otimes_\ZZ \RR$ dual to the fan $\Sigma$, the tautological ample line bundle on $X_\Sigma$ is denoted by $\CO_{X_\Sigma}(\Delta)$.

\subsubsection{Graphs}
The graphs we consider in this paper are finite connected graphs. They are allowed to have loops and multiple edges. For a given graph $\Gamma$, the sets of vertices and edges of $\Gamma$ are denoted by $V(\Gamma)$ and $E(\Gamma)$. For $v\in V(\Gamma)$, $val(v)$ denotes the valency of $v$. $V_k(\Gamma)$ denotes the set of vertices of valency $k$. If $v,v'\in V(\Gamma)$ then $E_{vv'}(\Gamma)$ denotes the set of edges connecting $v$ and $v'$. Most graphs in the paper are topological graphs, i.e., CW complexes of dimension one consisting of: (i) a 0-dimensional cell for each vertex, and (ii) a 1-dimensional cell for each edge glued to the 0-dimensional cells corresponding to the boundary vertices of the edge.

\section{Tropicalization}\label{sec:trop}
The canonical tropicalization procedure for curves over non-Archimedean fields in arbitrary characteristic was developed in \cite{T10}. In this section, we remind the construction and summarize the facts needed for the proof of the main result.

\subsection{Tropical and parameterized tropical curves}
Several different definitions of (parameterized) tropical curves can be found in the literature. Below, we follow \cite{T10}, and give a version of the definitions that are most convenient for our purposes.

\begin{definition}\label{def:trcur}
A {\it tropical curve} is a topological graph $\Gamma$ equipped with a complete, possibly degenerate, inner metric and with the following structure (s1),(s2); and satisfying the following properties (p1),(p2),(p3):
\begin{enumerate}
    \item[(s1)] $\Gamma$ has two types of vertices: {\it finite vertices} and {\it infinite vertices},
    \item[(s2)] the set of infinite vertices is equipped with a complete order, and is denoted by $V^\infty(\Gamma)$, the set of finite vertices is just a set and is denoted by $V^f(\Gamma)$;
    \item[(p1)] $\Gamma$ has finitely many vertices and edges;
    \item[(p2)] any infinite vertex has valency one and is connected to a finite vertex by an edge, called {\it unbounded edge}. Other edges are called {\it bounded edges}. The set of bounded edges is denoted by $E^b(\Gamma)$, and of unbounded edges by $E^\infty(\Gamma)$;
    \item[(p3)] any bounded edge $e$ is isometric to a closed interval $\left[0,|e|\right]$, where $|e|\in\RR$ denotes the length of $e$, and any unbounded edge $e$ is isometric to $[0,\infty]$, where the isometry maps the infinite vertex to $\infty$. Hence $|e|=\infty$ if $e$ is unbounded, and the restriction of the metric to $\Gamma\setminus V^\infty(\Gamma)$ is non-degenerate.
\end{enumerate}
A {\it $\QQ$-tropical curve} is a tropical curve such that $|e|\in\QQ\cup\{\infty\}$ for any $e\in E(\Gamma)$. A tropical curve is called {\it irreducible} if the underlying graph $\Gamma$ is connected. The connected components of $\Gamma$ are called {\em irreducible components}. The {\it genus} of a tropical curve $\Gamma$ is defined by $g(\Gamma):= 1-\chi(\Gamma)=1-|V(\Gamma)|+|E(\Gamma)|$. If $\Gamma$ is irreducible then $g(\Gamma)=b_1(\Gamma)$. A tropical curve is called {\it stable} if all its finite vertices have valency at least $3$. An {\it isomorphism} of tropical curves is an isomorphism of metric graphs.
\end{definition}

\begin{remark}\label{rem:trcur}
Let us explain the algebra-geometric motivation for this version of the definition:
Let $(C,D)$ be a smooth curve with marked points defined over the field $\FF$. Let $(C_{R_\LL}, D_{R_\LL})$ be a nodal model of $(C,D)$. One can associate to it a tropical curve $\Gamma_{C_{R_\LL},D_{R_\LL}}$ in the following way: The set of finite vertices is the set of irreducible components of the reduction of $C_{R_\LL}$, and the set of infinite vertices is the set of marked points $D\equiv D_{R_\LL}$.
The set of edges connecting two finite vertices is defined to be the set of common nodes of the corresponding components. In particular, if a component $C_v$ is singular then each singular point of $C_v$ corresponds to a loop at the corresponding finite vertex $v$. Finally, if a marked point specializes to certain component then the corresponding vertices are connected by an unbounded edge. It remains to specify the lengths of the bounded edges of $\Gamma_{C_{R_\LL},D_{R_\LL}}$. For a bounded edge $e$, set $|e|:=\frac{r_e+1}{e_\LL}$ if $C_{R_\LL}$ has singularity of type $A_{r_e}$ at the corresponding node, i.e., \'etale locally it is given by an equation $xy=t_\LL^{r_e+1}$. Observe that the length $|e|$ is independent of $\LL$. Indeed, if $\LL\subset\LL'$, $(C_{R_{\LL'}},D_{R_{\LL'}})=(C_{R_{\LL}},D_{R_{\LL}})\times_{\Spec R_\LL}\Spec R_{\LL'}$, and $C_{R_\LL}$ has singularity of type $A_r$ at a node $p$ then $C_{R_{\LL'}}$ has singularity of type $A_{e_{\LL'/\LL}(r+1)-1}$; hence $\frac{r+1}{e_\LL}=\frac{e_{\LL'/\LL}(r+1)-1+1}{e_{\LL'}}$. If the pair $(C,D)$ is stable then it admits a distinguished model, namely {\em the stable model,} and the associated tropical curve is independent of the field extension $\LL$. Note that $C$ is irreducible if and only if $\Gamma_{C_{R_\LL},D_{R_\LL}}$ is so. Note also, that the genus of $C$ is equal to the genus of $\Gamma_{C_{R_\LL},D_{R_\LL}}$ plus the sum of the genera of the irreducible components of the reduction of $C_{R_\LL}$. In particular, $g(C)=g(\Gamma_{C_{R_\LL},D_{R_\LL}})$ if and only if $C$ is a Mumford curve.
\end{remark}
\begin{example}\label{ex:trline}
Let $\FF:=\CCC((t))$ be the field of Laurent power series, $R:=\CCC[[t]]$ be the ring of integers, $C\subset\PP^2(\FF)$ be the line given by the homogeneous equation $x+ty=z$, and $q_1=[1-t:1:1], q_2=[1:0:1], q_3=[-t:1:0], q_4=[0:1:t]$ be the marked points. Let us describe the stable model of $(C,D)$ in this case: the homogeneous equation $x+ty=z$ defines a nodal model of the {\em curve} $C$ over the ring of integers $R$, whose reduction has unique component $L$ (the line in the complex plane given by the equation $x=z$). Plainly, $q_3$ and $q_4$ specialize to the same point $[0:1:0]$. Hence this integral model is not even a nodal model of the {\em pair} $(C,D)$. To resolve this issue, one must blow it up at the point $[0:1:0]$. In local coordinates, the initial model was given by $\Spec\,\CCC[[t]][\frac{x}{y},\frac{z}{y}]/(\frac{x}{y}-\frac{z}{y}+t)\simeq \Spec\,\CCC[[t]][\frac{x}{y}]$, the marked points $q_3, q_4$ were given by $\frac{x}{y}=-t$ and $\frac{x}{y}=0$, $L$ was given by $t=0$, and the blow up we perform has center at $t=\frac{x}{y}=0$. Denote the exceptional divisor by $E$. Now, $q_1, q_2$ specialize to two distinct points of the proper transform of $L$ (which we denote again by $L$), and $q_3,q_4$ specialize to two distinct points of $E$. Furthermore, the specializations are distinct from the node of the reduction. Thus, we have constructed a nodal model, which is stable since each component of the reduction contains three special points. We can now describe the tropical curve associated to the stable model of the pair $(C,D)$: It has two finite vertices $v_L$ and $v_E$ corresponding to the components $L$ and $E$ of the reduction. The finite vertices are joined by a unique bounded edge $e$ of length one, since the stable model is defined over $\FF$, and the intersection point $E\cap L$ is a regular point of the stable model. Furthermore, there are four infinite vertices corresponding to the marked points. The infinite vertices corresponding to $q_1$ and $q_2$ are connected to $v_L$, and the rest are connected to $v_E$.
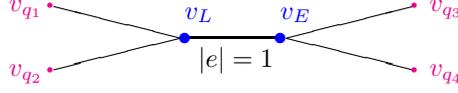
\begin{figure}[h]
\setlength{\unitlength}{0.1in} 
\centering 
\begin{picture}(20,6) 
\color{black}
\put(13.33,4){\line(4, 1){6.6}}
\put(13.33,4){\line(4, -1){6.6}}
\put(7.75,4){\line(-4, 1){6.6}}
\put(7.75,4){\line(-4, -1){6.6}}
\thicklines
\put(8.25,4){\line(1, 0){4.5}}
\put(8.7, 2.5){$|e|=1$}
\color{blue}
\put(8,4){\circle*{0.5}}
\put(13,4){\circle*{0.5}}
\put(8, 5){$v_L$}
\put(13, 5){$v_E$}
\color{magenta}
\put(0.95,2.3){\circle*{0.22}}
\put(20,2.3){\circle*{0.22}}
\put(0.95,5.7){\circle*{0.22}}
\put(20,5.7){\circle*{0.22}}
\put(20.8, 2){$v_{q_4}$}
\put(20.8, 5.4){$v_{q_3}$}
\put(-1.2, 2){$v_{q_2}$}
\put(-1.2, 5.4){$v_{q_1}$}
\end{picture}
\end{figure}
\end{example}

\begin{definition}\label{def:partrcur}
Let $N$ be a lattice. An {\it $N_\RR$-parameterized tropical curve} is a pair $(\Gamma, h_\Gamma)$, where $\Gamma$ is a tropical curve, and $h_\Gamma\colon V(\Gamma)\to N_\RR$ is a map such that
\begin{enumerate}
  \item $h_\Gamma(v)\in N$ for any infinite vertex $v\in V^\infty(\Gamma)$;
  \item $\frac{1}{|e|}(h_\Gamma(v)-h_\Gamma(v'))\in N$ for any bounded edge $e\in E_{vv'}(\Gamma)$;
  \item ({\it Balancing condition}) for any finite vertex $v$ the following holds: $$\sum_{v'\in V^f(\Gamma),\, e\in E_{vv'}(\Gamma)}\frac{1}{|e|}\left(h_\Gamma(v')-h_\Gamma(v)\right)+\sum_{v'\in V^\infty(\Gamma),\, e\in E_{vv'}(\Gamma)}h_\Gamma(v')=0.$$
\end{enumerate}
If $h_\Gamma(v)\in N_\QQ$ for all vertices $v$ then $\Gamma$ is called {\it $N_\QQ$-parameterized $\QQ$-tropical curve}.
\end{definition}
\begin{remark}
Usually, one defines a parameterized tropical curve to be a tropical curve $\Gamma$ together with a map $h\colon\Gamma\setminus V^\infty(\Gamma)\to N_\RR$ satisfying certain properties. Note, that after identifying the edges with straight intervals, a parameterized tropical curve in the sense of Definition~\ref{def:partrcur} defines $h$ as follows: $h$ is the unique continuous map that coincides with $h_\Gamma$ on the set of finite vertices, maps bounded edges $e\in E_{vv'}(\Gamma)$ linearly onto the intervals $[h_\Gamma(v),h_\Gamma(v')]$, and maps unbounded edges $e\in E_{vv'}(\Gamma)$, with $v'\in V^\infty(\Gamma)$, linearly onto the intervals $h_\Gamma(v)+\RR_+\cdot h_\Gamma(v')$.
\end{remark}

\begin{remark}\label{rem:partrcur}
Let us give some algebra-geometric motivation for this definition:
Let $(C,D)$ be a smooth curve with marked points defined over the field $\FF$, and let $f\colon C\setminus D\to T_N(\FF)$ be a morphism to an algebraic torus. Recall that if $(C,D)$ is stable then there is a tropical curve associated to it in a canonical way. Let $\Gamma$ be this tropical curve. We claim that $\Gamma$ admits a natural structure of an $N_\QQ$-parameterized $\QQ$-tropical curve. Indeed,let $v$ be a vertex. Then either $v$ corresponds to a component of the reduction, or it corresponds to a marked point. In both cases, the order of vanishing ${\rm ord}_v(f^*(x^m))$ is a linear function on $M$, hence an element of $N$. Set $h_\Gamma(v):=\frac{1}{e_\LL}{\rm ord}_v(f^*(x^\bullet))\in N_\QQ$ if $v$ is finite, and $h_\Gamma(v):={\rm ord}_v(f^*(x^\bullet))\in N$ if $v$ is infinite. Then $h_\Gamma$ is independent of the choice of $\LL$, and $(\Gamma, h_\Gamma)$ is an $N_\QQ$-parameterized $\QQ$-tropical curve by \cite[Lemma 2.11]{T10}.
\end{remark}
\begin{example}\label{ex:partrline}
Let $(C,D)$ be as in Example~\ref{ex:trline}, i.e., $C\subset\PP^2$ is the line given by the homogeneous equation $x+ty=z$ over the field $\FF:=\CCC((t))$, and the marked points are $q_1=[1-t:1:1], q_2=[1:0:1], q_3=[-t:1:0], q_4=[0:1:t]$. As we have seen in Example~\ref{ex:trline}, the tropical curve $\Gamma$ corresponding to the stable model of the pair $(C,D)$ has two finite vertices $v_L$ and $v_E$ joined by an edge of length one, and four infinite vertices. The infinite vertices corresponding to $q_1$ and $q_2$ are connected to $v_L$, and the rest are connected to $v_E$. Let $T:=\GG_m^2(\FF)\subset \PP^2(\FF)$ be the standard torus. Then we have a natural embedding $f\colon C\setminus D\to T$, and hence $\Gamma$ admits a natural structure of a parameterized tropical curve. Let us describe the corresponding function $h_\Gamma$:
The function $\frac{x}{z}$ is invertible at the generic point of $L$, and at the generic point of $E$, and at $q_1$, and at $q_2$. It vanishes to order one at $q_4$, and it has a simple pole at $q_3$. Similarly, the function $\frac{y}{z}$ is invertible at the generic point of $L$, and at $q_1$, and at $q_4$. It vanishes to order one at $q_2$, and it has a simple pole at the generic point of $E$ and at $q_3$. Thus, $h_\Gamma(v_L)=(0,0)$, $h_\Gamma(v_E)=(0,-1)$, $h_\Gamma(v_{q_1})=(0,0)$, $h_\Gamma(v_{q_2})=(0,1)$, $h_\Gamma(v_{q_3})=(-1,-1)$, $h_\Gamma(v_{q_4})=(1,0)$. In this case, it is very easy to verify the balancing condition directly: at $v_E$ we have $[(0,0)-(0,-1)]+(1,0)+(-1,-1)=(0,0)$, and at $v_L$ we have $[(0,-1)-(0,0)]+(0,1)=(0,0)$. Finally, let us describe the corresponding map $h\colon\Gamma\setminus V^\infty(\Gamma)\to \RR^2$: it contracts the unbounded edge containing $v_{q_1}$ to the point $(0,0)=h_\Gamma(v_L)$, maps the bounded edge to the straight interval joining $(0,0)$ and $(0,-1)$, and maps the unbounded edges connected to $q_2$, $q_3$, $q_4$ to the rays $(0,0)+\RR_+\cdot (0,1)$, $(0,-1)+\RR_+\cdot (-1,-1)$, $(0,-1)+\RR_+\cdot (1,0)$.
\end{example}
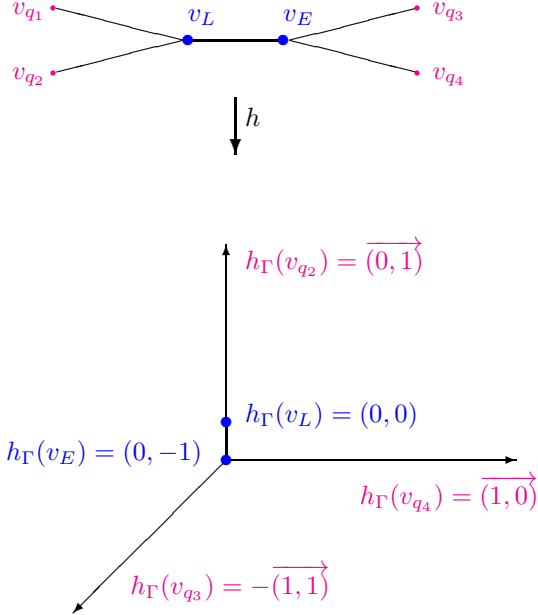
\begin{figure}[h]
\setlength{\unitlength}{0.1in} 
\centering 
\begin{picture}(20,36) 
\color{black}
\put(13.33,29){\line(4, 1){6.6}}
\put(13.33,29){\line(4, -1){6.6}}
\put(7.75,29){\line(-4, 1){6.6}}
\put(7.75,29){\line(-4, -1){6.6}}
\put(10.2,7){\vector(1, 0){15}}
\put(10,7){\vector(-1, -1){8}}
\put(10,9.3){\vector(0, 1){9}}
\thicklines
\put(8.25,29){\line(1, 0){4.5}}
\put(10.5,26){\vector(0, -1){3}}
\put(10,7.3){\line(0, 1){1.5}}
\put(11, 24.5){$h$}
\color{blue}
\put(8,29){\circle*{0.5}}
\put(13,29){\circle*{0.5}}
\put(10,7){\circle*{0.5}}
\put(10,9){\circle*{0.5}}
\put(8, 30){$v_L$}
\put(13, 30){$v_E$}
\put(11, 9){$h_\Gamma(v_L)=(0,0)$}
\put(-1.5, 7){$h_\Gamma(v_E)=(0,-1)$}
\color{magenta}
\put(0.95,27.3){\circle*{0.22}}
\put(20,27.3){\circle*{0.22}}
\put(0.95,30.7){\circle*{0.22}}
\put(20,30.7){\circle*{0.22}}
\put(20.8, 27){$v_{q_4}$}
\put(20.8, 30.4){$v_{q_3}$}
\put(-1.2, 27){$v_{q_2}$}
\put(-1.2, 30.4){$v_{q_1}$}
\put(11, 17){$h_\Gamma(v_{q_2})=\overrightarrow{(0,1)}$}
\put(17, 4.7){$h_\Gamma(v_{q_4})=\overrightarrow{(1,0)}$}
\put(5, 0){$h_\Gamma(v_{q_3})=-\overrightarrow{(1,1)}$}
\end{picture}
\end{figure}
\begin{remark}
It is easy to check that in Example~\ref{ex:partrline}, $h(\Gamma\setminus V^\infty(\Gamma))\cap\QQ^2$ coincides with the non-Archimedean amoeba $$\CA(C):=\left\{\left(\nu\left(\frac{x}{z}(p)\right), \nu\left(\frac{y}{z}(p)\right)\right)|p\in C(\overline{\FF})\cap \GG_m^2(\overline{\FF})\right\}.$$ One can prove that a similar statement holds true for any curve in a toric variety.
\end{remark}

\begin{definition}\label{def:proppartropcur} Let $(\Gamma, h_\Gamma)$ be an $N_\RR$-parameterized tropical curve, $v\in V^f(\Gamma)$ be a finite vertex, $e\in E_{vv'}(\Gamma)$ be an edge, and $v''\in V^\infty(\Gamma)$ be an infinite vertex.
\begin{enumerate}
  \item The {\it multiplicity} of $e$ is the integral length of $\frac{1}{|e|}(h_\Gamma(v)-h_\Gamma(v'))$ if $e$ is bounded, and is the integral length of $h_\Gamma(v')$ if $e$ is unbounded. The multiplicity of $e$ is denoted by $l(e)$.
  \item The {\it multiplicity} of $v''$ is the integral length of $h_\Gamma(v'')$.
  \item The {\it slope} of $e$ is the subspace $\RR\cdot(h_\Gamma(v)-h_\Gamma(v'))\subseteq N_\RR$ if $e$ is bounded, and is $\RR\cdot h_\Gamma(v')\subseteq N_\RR$ if $e$ is unbounded. The slope of $e$ is denoted by $N_{\RR,e}$, and the lattice $N\cap N_{\RR,e}$ is denoted by $N_e$. If $N_{\RR,e}\ne 0$ then $N_e$ and $N_{\RR,e}$ have a generator $n_e=\frac{1}{l(e)|e|}(h_\Gamma(v)-h_\Gamma(v'))$ if $e$ is bounded, and $n_e=\frac{1}{l(e)}h_\Gamma(v')$ if $e$ is unbounded. In the second case, it is a {\it distinguished} generator, while in the first case, it is defined only up-to a sign. However, if an {\em orientation} of the bounded edge is given then the generator is also distinguished.
  \item The {\it degree} of $\Gamma$ is the collection of pairs $(n_k,d_k)$, where $\{n_1,\dotsc ,n_s\}$ is the set of non-zero distinguished generators of slopes of unbounded edges, and $d_k=\sum_{e\in E^\infty(\Gamma), n_e=n_k}l(e)$. The degree is denoted by $\deg(\Gamma)$.
  \item By a {\em combinatorial type} of $(\Gamma, h_\Gamma)$ we mean the isomorphism class of the underlying graph $\Gamma$ equipped with the sublattice $l(e)N_e$ for all $e\in E(\Gamma)$.
\end{enumerate}
\end{definition}
\begin{remark}
Balancing condition implies $\sum_{(n,d)\in\deg(\Gamma)}dn=0$.
\end{remark}
\begin{example}
In Example~\ref{ex:partrline}, the multiplicities of all edges are one, $N_e=\ZZ\cdot(0,1)$,
$n_{e_1}=\overrightarrow{(0,0)}$, $n_{e_2}=\overrightarrow{(0,1)}$, $n_{e_3}=-\overrightarrow{(1,1)}$, and $n_{e_4}=\overrightarrow{(1,0)}$, where $e_i$ denotes the unbounded edge connected to the infinite vertex $v_{q_i}$ for any $i$. Finally, $\deg(\Gamma)=\{(n_{e_2},1), (n_{e_3},1), (n_{e_4},1)\}$.
\end{example}

The following combinatorial lemma is a version of Mikhalkin's \cite[Proposition~4.13]{Mikh05} and Nishinou-Siebert's \cite[Proposition 2.1]{NS06} that we will need in the proof of Theorem~\ref{thm:genZar} (1) in the toric setting. The lemma follows from either of these propositions, but since our definitions are slightly different, it may be not obvious. Thus, we include a proof for the convenience of the reader.
\begin{lemma}\label{lem:FinNumTyp}
For $r\in \NN$, the number of combinatorial types of $N_\RR$-parameterized stable tropical curves $(\Gamma, h_\Gamma)$ of given degree and genus, for which $|V^\infty(\Gamma)|<r$, is finite.
\end{lemma}
\begin{proof}
Fix a basis $\{m_i\}\subset M$. Then any $N_\RR$-parameterized tropical curve $(\Gamma, h_\Gamma)$ defines  $\RR$-parameterized tropical curves $(\Gamma, h^i_\Gamma):=(\Gamma, m_i\circ h_\Gamma)$. Fix an orientation of the bounded edges of $\Gamma$, and let $l^i(e)$ and $n_e^i$ be the multiplicity and the distinguished generator of the slope of $e$ in $(\Gamma, h^i_\Gamma)$. Then  $l(e)n_e=\sum_i l^i(e)n_e^i$. Thus, it is sufficient to prove the lemma for $N=\ZZ$. The degree in this case is $\deg=\{(\overrightarrow{1}, d),(\overrightarrow{-1}, d)\}$, which we denote simply by $d$; and the combinatorial type of $(\Gamma, h_\Gamma)$ is completely determined by the function $l\colon E(\Gamma)\to \ZZ_{\ge 0}$.

Let $(\Gamma, h_\Gamma)$ be an $\RR$-parameterized stable tropical curve of degree $d$ and genus $g$. Then $\sum_{e\in E^\infty(\Gamma)}l(e)=2d$. Fix an orientation of the edges of $\Gamma$, and a complete ordering on $V^f(\Gamma)$, both compatible with the standard orientation of $\RR$ via $h_\Gamma$. Then the balancing condition at a finite vertex $v$ is equivalent to the following: $\sum_{e\in E(\Gamma)}\epsilon(e,v)l(e)=0$, where $\epsilon(e,v)=-1$ if $v$ is the initial point of $e$, $\epsilon(e,v)=1$ if $v$ is the target of $e$, and $\epsilon(e,v)=0$ otherwise. Let us associate to $(\Gamma, h_\Gamma)$ the following datum: the isomorphism class of the underlying graph $\Gamma$, the orientation of the edges, the ordering on $V^f(\Gamma)$, and the function $l^\infty:=l|_{E^\infty(\Gamma)}$.

We claim that the set of such data associated to $\RR$-parameterized stable tropical curves of degree $d$ and genus $g$ is finite. Indeed, it is sufficient to show that the number of isomorphism classes of the underlying graphs $\Gamma$ is finite. But $|E(\Gamma)|=|V(\Gamma)|+g-1< 2r+3g-3$, since $|E(\Gamma)|-|V(\Gamma)|=g-1$ and $2|E(\Gamma)|=\sum_{v\in V(\Gamma)}val(v)>3|V(\Gamma)|-2r$, which implies the claim.

It remains to show that the number of the combinatorial types of $(\Gamma, h_\Gamma)$ corresponding to a given datum is finite; but the latter is plainly true by induction on the linearly ordered set $V^f(\Gamma)$, since the balancing condition allows only finitely many possibilities for $l(e)$ for any edge $e$ with initial vertex $v$ if the values $l(e)$ are given for all edges $e$ with initial vertices $u<v$.
\end{proof}

\subsection{The parameterized tropical curve corresponding to a reduced curve in a toric variety defined over $\FF$}
Let $X_\Sigma$ be a toric variety. Set $X:=X_\Sigma(\FF)$ and $T:=T_N(\FF)$. Let $x_1,\dotsc,x_k\in T$ be points, $C\subset X$ be a reduced irreducible curve containing $\{x_1,\dotsc,x_k\}$ in its smooth locus, and $f\colon \widetilde{C}\to X$ be the natural map from its normalization to $X$. Assume that $C$ does not intersect the union of orbits of codimension greater than one. Set $D:=f^{-1}(\{x_1,\dotsc,x_k\}\cup \partial X)$, and fix a linear ordering $D=\{q_1,\dotsc, q_r\}$ such that $x_i=f(q_i)$ for all $1\le i\le k$. Let $(C_{R_\LL}, D_{R_\LL})$ be the stable model of $(\widetilde{C}, D)$.

We define the corresponding $N_\QQ$-parameterized $\QQ$-tropical curve $(\Gamma, h_\Gamma)$ as follows:
$\Gamma$ is the tropical of the stable model of the pair $(\widetilde{C}, D)$ as constructed in Remark~\ref{rem:trcur}, and $(\Gamma, h_\Gamma)$ is the corresponding parameterized tropical curve as constructed in Remark~\ref{rem:partrcur}. The following claim follows immediately from the definitions:

\begin{claim}\label{cl:mult} The multiplicity of $q_i$ in $f^*(\partial X)$ is equal to the multiplicity of $v_{q_i}$. If $f(q_i)$ belongs to the big orbit of the divisor $E_\rho\subseteq \partial X$ then $h_\Gamma(v_{q_i})\in \rho$.
\end{claim}

Let $q_i\in D$ be a point, and $v_i$ be the unique finite vertex connected to $v_{q_i}$. Then $t_\LL^{-e_\LL h_\Gamma(v_i)(m)}f^*(x^m)$ does not vanish along the closure of $q_i$ in $D_{R_\LL}$ for all $m\in h_\Gamma(v_{q_i})^\perp\subset M$. Hence the following equality holds for all $m\in h_\Gamma(v_{q_i})^\perp\subset M$:
\begin{equation}
h_\Gamma(v_i)(m)=\nu(f^*(x^m)(q_i))=\nu(x^m(f(q_i)))
\end{equation}
From now on we assume that $f(q_i)$ belongs to the big orbits of $\partial X$ for $i>k$.
\begin{notation}\label{not:AL}
For $i\le k$, denote by $\CA_i\in N_\QQ$ the unique element mapping $m$ to $\nu(x^m(x_i))$. For $i>k$, set $L_{f(q_i)}:=\{n\in N_\QQ\, |\, n(m)=\nu(x^m(f(q_i))),\; \forall m\in \rho^\perp \}$, where $\rho\in\Sigma$ is such that $f(q_i)\in E_\rho$.
\end{notation}
\begin{corollary}\label{cor:constr}
$h_\Gamma(v_i)=\CA_i$ if $i\le k$, and $h_\Gamma(v_i)\in L_{f(q_i)}$ otherwise.
\end{corollary}
\begin{example}
In Example~\ref{ex:partrline}, $k=1$, $v_1=v_2=v_L$, and $v_3=v_4=v_E$. Plainly, $\CA_1=(0,0)$ since $\nu(\frac{x}{z}(q_1))=\nu(1-t)=0=\nu(1)=\nu(\frac{y}{z}(q_1))$. Thus, $\CA_1=h_\Gamma(v_1)$ as expected. The fan of the projective plane contains three rays $\rho_2:=\RR_+\cdot(0,1)$, $\rho_3:=\RR_+\cdot(-1,-1)$, and $\rho_4:=\RR_+\cdot(1,0)$, and for each $i>1$ the point $q_i$ belongs to the big orbit of the divisor $E_{\rho_i}$. Thus, by definition, we have $L_{f(q_2)}=\{(0,s)|s\in\RR\}$, $L_{f(q_3)}=\{(s+1,s)|s\in\RR\}$, $L_{f(q_4)}=\{(s,-1)|s\in\RR\}$; and $h_\Gamma(v_i)\in L_{f(q_i)}$ for all $i>1$ as expected.
\end{example}

\subsection{Deformations of parameterized tropical curves}
\begin{definition}\label{def:deforpartrcur} Let $(\Gamma, h_\Gamma)$ be an $N_\RR$-parameterized tropical curve. By a {\it deformation} of $(\Gamma, h_\Gamma)$ we mean a germ of a continuous family $\left\{\left(\Gamma_s,h_{\Gamma_s}\right)\right\}_{s\in(\RR, 0)}$ of $N_\RR$-parameterized tropical curves such that $\left(\Gamma_0, h_{\Gamma_0}\right)=(\Gamma, h_\Gamma)$, and the combinatorial type of the underlying graph of $\Gamma_s$ is independent of $s$.
\end{definition}
Any deformation of $(\Gamma, h_\Gamma)$ induces a deformation of the underlying graph, which can be canonically trivialized. Hence we may consider only deformations of $\Gamma$ inducing the trivial deformation of the underlying graph.
The multiplicities and the slopes of the edges are preserved by deformations since the lattice $N\subset N_\RR$ is discrete. This also shows that the lengths $|e|_s$ of bounded edges $e\in E_{vv'}$ with non-trivial slopes are uniquely defined by the values of $h_{\Gamma_s}(v)$ and $h_{\Gamma_s}(v')$, since the integral length $l(e)$ of $\frac{1}{|e|_s}\left(h_{\Gamma_s}(v)-h_{\Gamma_s}(v')\right)\in N$ is independent of $s$.

Fix an orientation of the bounded edges of $\Gamma$, and consider the linear map
$$\oplus_{v\in V^f(\Gamma)}N_\RR\to\oplus_{e\in E^b(\Gamma)}(N/N_e)_\RR$$
given by $x_v\mapsto\sum_{e\in E^b(\Gamma)}(\epsilon(e,v)x_v \bmod (N_e)_\RR)$, where $\epsilon(e,v)=-1$ if $v$ is the initial point of $e$, $\epsilon(e,v)=1$ if $v$ is the target of $e$, and $\epsilon(e,v)=0$ otherwise. Denote its kernel by $\EE_\RR^1(\Gamma)$. Then the {\it universal} deformation of $\Gamma$, i.e., the space of deformations up to isomorphism, can be identified naturally with the germ at the identity of the group $\EE_\RR^1(\Gamma)\times\RR_{\ge 0}^{c(\Gamma)}$, where $c(\Gamma)$ denotes the number of bounded edges of $\Gamma$ with trivial slope: To a deformation $\left\{\left(\Gamma_s,h_{\Gamma_s}\right)\right\}_{s\in\RR}$ one associates the collection  $$\left[\left(h_{\Gamma_s}(v)-h_\Gamma(v)\right)_{v\in V^f(\Gamma)}, \left(\frac{|e|_s}{|e|_0}\right)_{e, N_e=0}\right]\in \EE_\RR^1(\Gamma)\times\RR_{\ge 0}^{c(\Gamma)}\, .$$
\begin{remark}
Deformations of $N_\QQ$-parameterized $\QQ$-tropical curves are controlled by $\EE_\QQ^1(\Gamma)\times\QQ_{\ge 0}^{c(\Gamma)}$ defined similarly.
\end{remark}
\begin{example}
In Example~\ref{ex:partrline}, $c(\Gamma)=0$. Fix the orientation of the bounded edge $e$ such that $v_E$ is the initial point. Then the linear map $$\RR^2\oplus\RR^2=\oplus_{v\in V^f(\Gamma)}N_\RR\to \oplus_{e\in E^b(\Gamma)}(N/N_e)_\RR=\RR^2/\RR\cdot (0,1)=\RR$$ is given by $((a,b),(c,d))\mapsto c-a$. Hence the universal deformation of $(\Gamma,h_\Gamma)$ is the germ of the group $E^1_\RR(\Gamma)=\{((a,b),(a,d))|a,b,d\in\RR\}\simeq \RR^3$ at the identity.
\end{example}

Let us now explain the idea of the proof of Theorem~\ref{thm:genZar} (1) for toric surface $S=X_\Sigma(\kk)$ and boundary divisor $E=\partial X_\Sigma(\kk)=\cup_\rho E_\rho(\kk)$. In this case $K_S+E=0$. Assume for simplicity that $|\alpha|=0$. If $V$ has dimension $k$ then there exists a curve of geometric genus $g$ in the linear system $\CL$ passing through $k$ general points $x_1,\dotsc, x_k$ in the torus, and we may choose these points such that the points $\CA_1,\cdots,\CA_k$ of Notation~\ref{not:AL} are in general position in the plane $N_\QQ$. Thus, the corresponding parameterized tropical curve $(\Gamma,h_\Gamma)$ ``passes" through the points $\CA_1,\cdots,\CA_k$, i.e., the first assertion of Corollary~\ref{cor:constr} is satisfied. However, by Mikhalkin's \cite[Proposition 2.23]{Mikh05}, the dimension of the universal deformation space of such tropical curve is at most $|\beta|+p_g(C)-1$, and hence $(\Gamma,h_\Gamma)$ may ``pass" through at most $|\beta|+p_g(C)-1$ points in general position. Thus, $\dim(V)=k\le |\beta|+p_g(C)-1=-C.(K_S+E)+|\beta|+p_g(C)-1$.

\section{The proof of the main theorem}\label{sec:mainthmproof}
In this section we prove Theorem~\ref{thm:genZar} (1) for a toric surface $S=X_\Sigma(\kk)$ and boundary divisor $E=\partial X_\Sigma(\kk)=\cup_\rho E_\rho(\kk)$.

Assume to the contrary that $\dim(V)\ge -C.(K_S+E)+|\beta|+p_g(C)=|\beta|+p_g(C)$. Set $k:=\dim(V)$, $R:=\kk[[t]]$, $\FF:=Frac(R)$, and consider $S(\overline{\FF}),E(\overline{\FF})$, and $V(\overline{\FF})$. Then the curve corresponding to a general point of $V(\overline{\FF})$ satisfies the assumptions of Theorem~\ref{thm:genZar}. To simplify the notation we will omit $\overline{\FF}$ below.

For any irreducible component $E'=E_\rho\subset E$, the points $p_{i,j}^{E'}$ belong to the maximal orbit of $E'$, and we may assume that the collection of lines $L_{p_{i,j}^{E'}}$ of Notation~\ref{not:AL} is general in the set of lines with the slope $\rho$. Let $x_i\in T_N(\overline{\FF})$, $1\le i\le k$, be general points such that the collection $\{\CA_i\in N_\QQ\}$ of Notation~\ref{not:AL} is general. Then there exists a point in $V$ such that the corresponding curve contains $x_1,\dotsc,x_k$ in its smooth locus, and without loss of generality we may assume that this curve is $C$. After replacing $(R,\FF)$ with a finite separable extension, we may assume that $C$, $x_i$, and $p_{i,j}^{E'}$ are defined over $\FF$.

Consider the natural map $f\colon \widetilde{C}\to S$ from the normalization of $C$ to $S$, and set $D:=f^{-1}(\{x_1,\dotsc,x_k\}\cup E)$. Set $r:=k+|\alpha|+|\beta|$, and fix a linear ordering $D=\{q_1,\dotsc, q_r\}$, such that the first $k$ points are mapped to $x_1,\dotsc,x_k$, and the next $|\alpha|$ points are mapped to $\{p_{i,j}^{E'}\}_{E',i,j}$. Let $(\Gamma, h_\Gamma)$ be the $N_\QQ$-parameterized $\QQ$-tropical curve corresponding to $C$. Then $g(\Gamma)\le p_g(C)$, and $\deg(\Gamma)=\{(n_\rho,d_\rho)\}$, where $d_\rho=C.E_\rho$ by Calim~\ref{cl:mult}. Furthermore, $h_\Gamma(v_i)=\CA_i$ for $i\le k$ and $h_\Gamma(v_i)\in L_{f(q_i)}$ for $k<i\le k+|\alpha|$ by Corollary~\ref{cor:constr}.

By Lemma \ref{lem:FinNumTyp}, there exist only finitely many combinatorial types of $N_\QQ$-parame\-terized $\QQ$-tropical curves of degree $\{(n_\rho,d_\rho)\}$ and genus at most $p_g(C)$. Thus, without loss of generality, we may assume that for sufficiently small deformations $\CA_i(s)\in N_\QQ$ and $L_i(s)||L_{f(q_i)}$, there exists a deformation $\left(\Gamma_s,h_{\Gamma_s}\right)$ such that $h_{\Gamma_s}(v_i)=\CA_i(s)$ for $i\le k$ and $h_{\Gamma_s}(v_i)\in L_i(s)$ for $k<i\le k+|\alpha|$. One could finish the proof by saying that we get a contradiction to \cite[Proposition 2.23]{Mikh05}. However, to make the presentation self-contained, we give a complete proof here. It is different from Mikhalkin's proof of \cite[Proposition 2.23]{Mikh05}, but is based on the idea of the proof of \cite[Proposition 4.19]{Mikh05}.

Fix an orientation of the bounded edges of $\Gamma$. Then the natural linear projection
\begin{equation}\label{eq:surj}
\EE_\QQ^1(\Gamma)\twoheadrightarrow \oplus_{i=1}^kN_\QQ\bigoplus\oplus_{i=k+1}^r(N/N_{e_i})_\QQ\,
\end{equation}
is surjective. Consider the graph $\Gamma'$ obtained from $\Gamma$ by removing the vertices $v_{q_i}$ and $v_i$, and the unbounded edges $e_i=\overline{v_{q_i}v_i}$ for $i\le k$, and gluing one-valent finite vertices to the edges that have contained the finite vertices $v_i$. We note here that $(\Gamma', h_\Gamma|_{\Gamma'})$ is not a parameterized tropical curve since it does not satisfy the balancing condition at the new finite vertices. Then
$$-|\beta|>-1-|\beta|\ge -\chi(\Gamma)-k\ge -\chi(\Gamma')=-\sum\chi(\Gamma_j)=\sum \left(b_1(\Gamma_j)-1\right),$$ where $\Gamma_j\subset \Gamma'$ are the connected components. Thus, there exists $j$ such that $b_1(\Gamma_j)=0$, and $v_{q_i}\notin V(\Gamma_j)$ for all $k+|\alpha|<i\le r=k+|\alpha|+|\beta|$. Without loss of generality we may assume that $j=1$. Since $b_1(\Gamma_1)=0$, it follows (e.g. by induction on $|V^f(\Gamma_1)|$) that the natural map $\oplus_{v\in V^f(\Gamma_1)}N_\QQ\to\oplus_{e\in E^b(\Gamma_1)}N_\QQ$ is surjective; hence $\oplus_{v\in V^f(\Gamma_1)}N_\QQ\to\oplus_{e\in E^b(\Gamma_1)}(N/N_e)_\QQ$ is also surjective. Thus,
$$\oplus_{v\in V^f(\Gamma_1)}N_\QQ\to\oplus_{e\in E^b(\Gamma_1)}(N/N_e)_\QQ\bigoplus\oplus_{v\in V_1^f(\Gamma_1)}N_\QQ\bigoplus\oplus_{v\in E^\infty(\Gamma_1)}(N/N_e)_\QQ$$
is so by the surjectivity of \eqref{eq:surj}, where $V_1$ (resp. $V_1^f$) denotes the set of (resp. finite) vertices of valency one. Then $2|V^f(\Gamma_1)|\ge |E^b(\Gamma_1)|+2|V_1^f(\Gamma_1)|+|E^\infty(\Gamma_1)|$, or equivalently, $$2|V(\Gamma_1)|\ge |E(\Gamma_1)|+2|V_1(\Gamma_1)|.$$ On the other hand, $\Gamma$ is stable, hence $val(v)\ge 3$ for any $v\in V(\Gamma_1)\setminus V_1(\Gamma_1)$. Thus, $$2|E(\Gamma_1)|\ge 3|V(\Gamma_1)|-2|V_1(\Gamma_1)|.$$ The two inequalities imply $|E(\Gamma_1)|\ge |V(\Gamma_1)|$. However, $|V(\Gamma_1)|-|E(\Gamma_1)|=1$, since $b_1(\Gamma_1)=0$. We got a contradiction, and the proof is now complete.

\begin{remark} Similar argument can be used to obtain another proof of \cite[Proposition 2.23]{Mikh05}. Indeed, if $k=\dim(V)=|\beta|+p_g(C)-1$ then, to avoid the contradiction, one must have equalities in all the inequalities above. Thus, $b_1(\Gamma_j)\le 1$ for any $j$. If $b_1(\Gamma_j)=0$ then $V(\Gamma_j)$ contains unique $v_{q_i}$ with $k+|\alpha|<i\le r=k+|\alpha|+|\beta|$, no edge of $\Gamma_j$ has trivial slope, and $\Gamma_j$ has only vertices of valency $3$ and $1$. Hence, if there exists $j$ such that $b_1(\Gamma_j)=1$ then $v_{q_i}\notin V(\Gamma_j)$ for all $k+|\alpha|<i\le r=k+|\alpha|+|\beta|$, and one easily gets a contradiction similar to the contradiction we got in the proof. Thus, $b_1(\Gamma_j)=0$ for all $j$, $\Gamma$ is trivalent, and no edge of $\Gamma$ has trivial slope.
\end{remark}
\section{Examples and counterexamples}\label{sec:exandcex}
{\em Throughout this section, $\kk$ denotes an algebraically closed ground field of positive characteristic $p$, $r\in\NN$, and $q=p^r$.} Let us fix the notation. Set $N:=\ZZ^2$, and
$$n_{k,1}:=  \left(
        \begin{smallmatrix}
          0 \\
          1 \\
        \end{smallmatrix}
      \right),
n_{k,2}:=  \left(
        \begin{smallmatrix}
          k \\
          1 \\
        \end{smallmatrix}
      \right),
n_{k,3}:=  \left(
        \begin{smallmatrix}
          -k \\
          -2 \\
        \end{smallmatrix}
      \right),
m_{k,1}:=\left(
        \begin{smallmatrix}
          1 \\
          0 \\
        \end{smallmatrix}
      \right),
m_{k,2}:=\left(
        \begin{smallmatrix}
          -1 \\
          k \\
        \end{smallmatrix}
      \right),
m_{k,3}:=\left(
        \begin{smallmatrix}
          0 \\
          0 \\
        \end{smallmatrix}
      \right);$$
$$n'_{k,1}:=  \left(
        \begin{smallmatrix}
          0 \\
          1 \\
        \end{smallmatrix}
      \right),
n'_{k,2}:=  \left(
        \begin{smallmatrix}
          k \\
          1 \\
        \end{smallmatrix}
      \right),
n'_{k,3}:=  \left(
        \begin{smallmatrix}
          0 \\
          -1 \\
        \end{smallmatrix}
      \right),
n'_{k,4}:=  \left(
        \begin{smallmatrix}
          -k \\
          -1 \\
        \end{smallmatrix}
      \right),$$
$$m'_{k,1}:=\left(
        \begin{smallmatrix}
          0 \\
          0 \\
        \end{smallmatrix}
      \right),
m'_{k,2}:=\left(
        \begin{smallmatrix}
          1 \\
          0 \\
        \end{smallmatrix}
      \right),
m'_{k,3}:=\left(
        \begin{smallmatrix}
          0 \\
          k \\
        \end{smallmatrix}
      \right),
m'_{k,4}:=\left(
        \begin{smallmatrix}
          -1 \\
          k \\
        \end{smallmatrix}
      \right).$$
Let $\Delta_k\subset M_\RR$ denote the triangle with vertices $m_{k,i}$, and $\Delta'_k\subset M_\RR$ denote the parallelogram with vertices $m'_{k,i}$. Let $\Sigma_k$ and $\Sigma'_k$ be the fans in $N_\RR$ dual to $\Delta_k$ and $\Delta'_k$, i.e., the complete fans generated by the rays $\rho_{k,i}=\RR_+n_{k,i}$ and $\rho'_{k,i}=\RR_+n'_{k,i}$ respectively. Set $S_k:=X_{\Sigma_k}(\kk)$, $S'_k:=X_{\Sigma'_k}(\kk)$, $\CL_k:=\CO_{S_k}(\Delta_k)$, and $\CL'_k:=\CO_{S'_k}(\Delta'_k)$. To simplify the notation, we will write $E_{\rho_i}$ and $E_{\rho'_i}$ instead of $E_{\rho_i}(\kk)$ and $E_{\rho'_i}(\kk)$ when refereeing to the components of the boundary divisors. We will also omit the subindex $k$ everywhere if $k$ is given and no confusion is possible.

\subsection{Counterexamples to Theorem~\ref{thm:genZar} (2) over $\kk$}
Two series of counterexamples are constructed in Theorems \ref{thm:char>2} and \ref{thm:chararb}. One can prove each of them either by a straight-forward computation, or by developing a more general approach. To demonstrate both, we prove Theorem \ref{thm:char>2} by a computation, and give a conceptual proof of Theorem \ref{thm:chararb}.

\begin{theorem}\label{thm:char>2}
Assume that $p>2$. Let $V^{irr}(S_q,\CL_q,0)$ be the Severi variety parameterizing irreducible rational curves in $|\CL_q|$ that do not contain the zero-dimensional orbits in $S_q$, and $C$ be a curve corresponding to a closed point of $V^{irr}(S_q,\CL_q,0)$. Then $C$ has a unique singular point, and the type of the singularity is $A_{q-1}$. If $C'$ corresponds to a general closed point of $V^{irr}(S_q,\CL_q,0)$ then $C'$ intersects $C$ at a unique point, and the intersection index at this point is $q$. Furthermore, $V^{irr}(S_q,\CL_q,0)$ is irreducible and has expected dimension $-C.K_{S_q}-1$.
\end{theorem}
\begin{proof}
Since $q$ is odd, the integral lengths of the sides of $\Delta_q$ are one, and hence $C.E_{\rho_i}=1$ for all $i$. Thus, there exists unique isomorphism $\PP^1\to \widetilde{C}$, where $\widetilde{C}$ denotes the normalization of $C$, taking $0,1$, and $\infty$ to the preimages of $C\cap E_{\rho_1}$, $C\cap E_{\rho_2}$, and $C\cap E_{\rho_3}$. Let $t$ be the coordinate on $\AAA^1=\PP^1\setminus\{\infty\}$. Then the map $f\colon \PP^1\to\widetilde{C}\to S_q$ is given by
\begin{equation}\label{eq:pullback}
f^*(x^m)=\chi(m)t^{(n_1,m)}(t-1)^{(n_2,m)},
\end{equation}
where $\chi\colon M\to \kk^*$ is a multiplicative character. Vice versa, for any multiplicative character $\chi\colon M\to \kk^*$, \eqref{eq:pullback} defines a map from $\PP^1$ to $S_q$ whose image corresponds to a closed point of $V^{irr}(S_q,\CL_q,0)$. We constructed a natural isomorphism $$\iota\colon V^{irr}(S_q,\CL_q,0)\to T_N(\kk).$$ Hence $V^{irr}(S_q,\CL_q,0)$ is irreducible, and has expected dimension $2=-C.K_{S_q}-1$.

We claim that the differential of $f$ vanishes at $t=\frac{1}{2}$. Indeed, the log derivatives of $f^*(x^m)$ are given by $\frac{(n_1,m)}{t}+\frac{(n_2,m)}{t-1}$; hence vanish at $t=\frac{1}{2}$. Let $e_1,e_2\in M$ be the standard basis, and $a:=x^{e_1}-x^{e_1}(f(\frac{1}{2})), b:=x^{e_2}-x^{e_2}(f(\frac{1}{2}))$ be the local coordinates at $f(\frac{1}{2})$. Then $$f^*(a)=\chi(e_1)\left(t-\frac{1}{2}\right)^q,\, f^*(b)=\chi(e_2)\left(t-\frac{1}{2}\right)^2,$$ and  the image of $f$ satisfies the equation $$\chi(e_2)^qa^2=\chi(e_1)^2b^q.$$ Hence the singularity of $C$ at $f(\frac{1}{2})$ is of type $A_{q-1}$. Note that the arithmetic genus of $\CL_q$ is equal to the number of integral points in the interior of $\Delta_q$, which is equal to $\frac{q-1}{2}=\delta(A_{q-1})$. Thus, $C$ has a unique singular point.

Since $\CL_q^2=2Area(\Delta_q)=q$, it is sufficient to show that $C\cap C'$ is a point. Set  ${\chi':=\iota(C')}$. Then $f'\colon \PP^1\to\widetilde{C'}\to S_p$ is given by \eqref{eq:pullback} with $\chi$ replaced by $\chi'$. If $f(s)=f'(s')$ then $f'^*(m)(s')=f^*(m)(s)$ for all $m\in M$. In particular, $$\chi(m_1)(s^q-1)=\chi(m_1)(s-1)^q=\chi'(m_1)(s'-1)^q=\chi'(m_1)((s')^q-1)$$ and $\chi(m_2)s^q=\chi'(m_2)(s')^q.$ Hence the intersection point is unique, and is given by $$s=\sqrt[q]{\frac{\chi'(m_2)\chi(m_1)-\chi'(m_2)\chi'(m_1)}{\chi(m_1)\chi'(m_2)-\chi'(m_1)\chi(m_2)}}\hspace{0.2cm}\text{and}\hspace{0.2cm} s'=\sqrt[q]{\frac{\chi(m_2)\chi(m_1)-\chi(m_2)\chi'(m_1)}{\chi(m_1)\chi'(m_2)-\chi'(m_1)\chi(m_2)}}.$$
The proof is now complete.
\end{proof}

\begin{theorem}\label{thm:chararb}
Let $V^{irr}(S'_q,\CL'_q,0)$ be the Severi variety parameterizing irreducible rational curves in $|\CL'_q|$ that do not contain the zero-dimensional orbits in $S'_q$, and $C$ be a curve corresponding to a closed point of $V^{irr}(S'_q,\CL'_q,0)$. Then all singularities of $C$ are unibranch, and the number of singular points of $C$ is one if $p=2$, and is two if $p>2$. If $C'$ corresponds to a general closed point of $V^{irr}(S'_q,\CL'_q,0)$ then $C'$ intersects $C$ at two points, and the intersection index at these points is $q$. Furthermore, $V^{irr}(S'_q,\CL'_q,0)$ is irreducible and has expected dimension $-C.K_{S'_q}-1$.
\end{theorem}
\begin{proof}
Consider the sublattice $N'\subset N$ spanned by $\left(
        \begin{smallmatrix}
          q \\
          0 \\
        \end{smallmatrix}
      \right)$ and $\left(
        \begin{smallmatrix}
          0 \\
          1 \\
        \end{smallmatrix}
      \right)$. Then the toric surfaces corresponding to the fan $\Sigma'_2$ with respect to $N'$ is $\PP^1\times\PP^1$, $S'_q$ is the quotient of $\PP^1\times\PP^1$ by the action of $\mu_q$, and $\pi\colon\PP^1\times\PP^1\to S'_q$ is bijective. Furthermore, since the restrictions of $N$ and of $N'$ onto the rays of $\Sigma'_2$ coincide, the restriction of the action onto the one-dimensional torus orbits of $\PP^1\times \PP^1$ is free.

Set $D:=\pi^{-1}(C)=C\times_{S'_q}(\PP^1\times\PP^1)$. The integral length of the sides of $\Delta'_q$ is one, hence $C.E_{\rho'_i}=1$, and $D$ intersects each coordinate line at a unique point with multiplicity $q$. Hence $D\in |\CO_{\PP^1\times\PP^1}(q,q)|$. If $\eta\in C$ is the generic point then $\eta\times_{S'_q}(\PP^1\times\PP^1)=\Spec \left(\kk(\eta)[(x^{m'_2})^{1/q}]\right)$. Thus, the reduced curve $L:=D^{\rm red}$ is rational and belongs to the linear system $|\CO_{\PP^1\times\PP^1}(q_1,q_1)|$ for some $q_1=p^{r_1}$. Furthermore, it intersects each coordinate line at a unique point, and is unibranch at these points.

Fix an isomorphism from $\PP^1$ to the normalization of $L$, and consider the projections $\pi_1,\pi_2\colon \PP^1\to L\to \PP^1\times\PP^1\rightrightarrows \PP^1$. Then $|\pi_i^{-1}(0)|=|\pi_i^{-1}(\infty)|=1$. Hence $\pi_i=Fr_{q_1}\circ\phi_i=\phi_i\circ Fr_{q_1}$ for some $\phi_i\in Aut(\PP^1)$, where $Fr_{q_1}\colon\PP^1\to\PP^1$ denotes the Frobenius morphism. Thus, $L\in |\CO_{\PP^1\times\PP^1}(1,1)|$,  $q_1=1$, and $L$ is the normalization of $C$. We constructed a bijection between $V^{irr}(S'_q,\CL'_q,0)$ and an open subset of $|\CO_{\PP^1\times\PP^1}(1,1)|$, hence $V^{irr}(S'_q,\CL'_q,0)$ is irreducible, and has dimension $3$ as expected.

Since two general curves in $|\CO_{\PP^1\times\PP^1}(1,1)|$ intersect at two points, it follows that their images on $S'_q$ are tangent to each other at two points to order $q$. Since $\pi$ is bijective, all singularities of $C$ are unibranch. Consider the (non-cartesian) diagram
$$\xymatrix{
L\ar@{^{(}->}[d]_\iota\ar[r]^f\ar[dr]^{\pi\iota} & C\ar@{^{(}->}[d]\\
\PP^1\times\PP^1 \ar[r]_-{\pi} & S'_q
}$$
and the corresponding diagram of sheaves of log-differential forms
$$\xymatrix{
M'\otimes_\ZZ\CO_L\ar@{=}[r] & \iota^*\Omega_{\PP^1\times\PP^1}\left(\log\left(\partial(\PP^1\times\PP^1)\right)\right)\ar[r]^-{d\iota} & \Omega_L\left(\log\left(\iota^*\partial(\PP^1\times\PP^1)\right)\right)\ar@{=}[d]\\
M\otimes_\ZZ\CO_L\ar@{=}[r]\ar[u] & (\pi\iota)^*\Omega_{S'_q}\left(\log\left(\partial(S'_q)\right)\right)\ar[r]^-{d(\pi\iota)} & \Omega_L\left(\log\left((\pi\iota)^*\partial(S'_q)\right)\right)\\
}$$
which exists since $L$ and $C$ intersect $\partial(\PP^1\times\PP^1)$ and $\partial(S'_q)$ transversally. For the same reason there exists a natural exact sequence
$$\xymatrix{
0\ar[r] & \CO_L(-L)\ar[r] & M'\otimes_\ZZ\CO_L\ar[r]^-{d\iota} & \Omega_L\left(\log\left(\iota^*\partial(\PP^1\times\PP^1)\right)\right)\ar[r] & 0,\\
}$$
and $\Omega_L\left(\log\left((\pi\iota)^*\partial(S'_q)\right)\right)=\Omega_L\left(\log\left(\iota^*\partial(\PP^1\times\PP^1)\right)\right)\simeq \CO_L(2)$, since $L^2=2$.

Let $z\in L$ be a point, and set $\CF:={\rm Coker}\left(d(\pi\iota)\right)$. Then $f(z)\in C$ is singular if and only if $z\in {\rm supp}(\CF)$. Hence the number of singular points of $C$ is equal to $|{\rm supp}(\CF)|$. The image of $M\otimes_\ZZ\CO_L$ in $M'\otimes_\ZZ\CO_L$ is isomorphic to $\CO_L$, since $\ZZ^2=M\to M'=\ZZ^2$ is given by $\left(
                                          \begin{smallmatrix}
                                            q & 0 \\
                                            0 & 1 \\
                                          \end{smallmatrix}
                                        \right)$.
Hence $ht(\CF)=2$. Furthermore, since $x^{\left(
                                          \begin{smallmatrix}
                                             0 \\
                                             1 \\
                                          \end{smallmatrix}
                                        \right)}|_L$
has two simple zeroes and two simple poles, $|{\rm supp}(\CF)|$ is equal to the number of zeroes of the differential $d\left(\frac{t-\xi}{t(t-1)}\right)$ on $\PP^1$ for some $\xi\ne 0,1,\infty$, which is equal to one if $p=2$, and is equal to two if $p>2$.
\end{proof}

\begin{remark}
(1) One can describe the type of the singularities of $C\subset S'_q$ in terms of Greuel-Kr\"oning classification of simple singularities in positive characteristic \cite{GK90}. If $p>2$ then the singularities of $C$ are of type $A_{q-1}$, and if $p=2$ then the singularity is of type $A_{2q-2}^{q-2}$, i.e. in formal local coordinates given by the equation $x^2+xy^q+y^{2q-1}=0$ or, equivalently, by $x^2+xy^q+\lambda y^q=0$,  $\lambda\ne 0$. To see this, one can do a straightforward computation similar to the computation in the second paragraph of the proof of Theorem~\ref{thm:char>2}, which we leave to the reader.

(2) Observe that the weighted projective plane $S_q$ is the quotient $\PP^2/\mu_q$. One can check that the rational curves $C,C'\subset S_q$ are the images of lines $L,L'\subset \PP^2$ that do not contain the zero-dimensional orbits. This observation can be used to give another proof of Theorem~\ref{thm:char>2} similar to the proof of Theorem~\ref{thm:chararb}.

(3) Let $k$ be a natural number divisible by $p$, and suppose that $q$ is the maximal power of $p$ dividing $k$. Consider the toric surfaces $S_k$ and $S'_k$, and the Severi varieties $V^{irr}(S_k,\CL_k,0)$ and $V^{irr}(S'_k,\CL'_k,0)$. One can show that the general point of $V^{irr}(S'_k,\CL'_k,0)$ corresponds to a curve with $\frac{k}{q}-1$ points of type $A_{2q-1}$, and either two singular points of type $A_{q-1}$ if $p>2$, or one singular point of type $A_{2q-2}^{q-2}$ if $p=2$. Furthermore, the curves corresponding to two general points in $V^{irr}(S'_k,\CL'_k,0)$ intersect at $2\frac{k}{q}$ distinct points, and have contact of order $q$ at each of them. Similarly, if $k$ is odd then one can show that the general closed point of $V^{irr}(S_k,\CL_k,0)$ corresponds to a curve with one singular point of type $A_{q-1}$, and $\frac{1}{2}\left(\frac{k}{q}-1\right)$ points of type $A_{2q-1}$. Furthermore, the curves corresponding to two general closed points in $V^{irr}(S_k,\CL_k,0)$ intersect at $\frac{k}{q}$ distinct points, and have contact of order $q$ at each of them. We leave the details to the reader.

(4) Another interesting example in characteristic $p=2$ can be obtained if, for $k$ divisible by four, one considers the Severi variety parameterizing rational curves in the linear system $|\CL_k|$ having unique intersection with each one of the three one-dimensional orbits, and not containing the zero-dimensional orbits. I this case, the general curve $C$ has a unibranch singularity at its point of intersection with $E_{\rho_3}$, and $\frac{1}{2}\left(\frac{k}{q}-1\right)$ points of type $A_{2q-1}$. As before, two general curves intersect non-transversally, and the Severi variety is irreducible and has expected dimension.

(5) The curves corresponding to the closed points of $V^{irr}(S_3,\CL_3,0)$ can be identified naturally with the plane cubic curves having contact of order three to two coordinate axis at their points of intersection with the third axis. To see this, observe that the triangle $\Delta_3$ is contained in the triangle with vertices
$
\left(
  \begin{smallmatrix}
    -1 \\
    0 \\
  \end{smallmatrix}
\right),
\left(
    \begin{smallmatrix}
        -1 \\
        3 \\
    \end{smallmatrix}
\right),
\left(
    \begin{smallmatrix}
        2 \\
        0 \\
     \end{smallmatrix}
\right)$,
whose integral points correspond to a basis of the space of plane cubic curves. Plainly, the integral points of $\Delta_3$ correspond to the basis of the subspace mentioned above. Assume that $p=3$. Then any curve corresponding to a closed point of $V^{irr}(S_3,\CL_3,0)$ has a cusp in the maximal orbit and any two such curves have a contact of order three by Theorem~\ref{thm:char>2}. Hence the same is true for plane cubic curves having contact of order three to two coordinate axis at their points of intersection with the third axis. Finally, since $Aut(\PP^2)$ acts transitively on the set of lines, it follows that the general curve in the Severi variety parameterizing {\em plane} cubics having contact of order three at {\em unspecified} points to a pair of lines is cuspidal, and any two such cubics have contact of order three. Another way to see that such cubics are cuspidal was suggested to us by Joe Harris: if a plane cubic $C$ contacts two lines to order three then there is an element of order three in the Picard group $Pic(C)$. However, $\GG_m$ has no elements of order equal to the characteristic. Hence $Pic(C)=\GG_a$, and $C$ is cuspidal.
\end{remark}

\subsection{Examples of reducible Severi varieties on toric surfaces over $\kk$}

\begin{theorem}\label{thm:redSev}
Let $d\ge 2$, $\frac{q-1}{2}\le g_1\le \min\left\{\frac{2dq-2d-q-1}{2}, \frac{(d-1)(d-2)}{2}\right\}$, and $q-1\le g_2\le \min\left\{2dq-q-d-1, (d-1)^2\right\}$ be positive integers. Then the Severi variety $V^{irr}(S'_q,(\CL'_q)^d, g_2)$ is reducible, and $V^{irr}(S_q,\CL_q^d, g_1)$ is reducible if $p>2$.
\end{theorem}

The idea of the proof is to show that the Severi varieties contain components whose general closed points correspond to curves with different types of singularities. To prove the theorem we will need some preparations.

\subsubsection{Deformation theory}\label{subsec:def}
\begin{lemma}\label{lem:defcur}
Let $k\ge 0$ be an integer, $S$ be a (not necessarily complete) smooth rational surface, $C\subset S$ be a complete reduced curve, $x_1,\dotsc,x_k\in C$ be some of its nodes, $C'$ be the partial normalization of $C$ preserving the nodes  $x_1,\dotsc,x_k$, and $f\colon C'\to S$ be the natural morphism. Assume that  $-f(E').K_S-\sum_{x\in E'}{\rm ord}_x(df)>0$ for any irreducible component $E'\subseteq C'$, where ${\rm ord}_x(df)$ denotes the order of vanishing of $df$ at $x$. Then the deformation space of the pair $(C',f)$ is smooth, unobstructed, and has expected dimension $-C.K_S+p_a(C')-1$. In particular, $C$ deforms to a curve of geometric genus $p_a(C')$ smooth in the neighborhood of the points $x_1,\dotsc,x_k$. Furthermore, if $C'$ is connected then the deformed curve is irreducible.
\end{lemma}
\begin{proof}
Let $\CN_f=Coker\left(\Theta_{C'}\to f^*\Theta_S\right)$ be the normal sheaf of $f$. Recall that the deformations of nodal curves are unobstructed, and any deformation of the nodes lifts to a global deformation of the curve. Furthermore, since $S$ is smooth, by the infinitesimal lifting property, we have a natural exact sequence
\begin{equation}\label{eq:seqdef}
0\to H^0(C',\CN_f)\to Def^1(C',f)\to H^0\left(C', \CE xt^1_{\CO_{C'}}(\Omega_{C'}, \CO_{C'})\right)\to 0\, ,
\end{equation}
and if $H^1(C',\CN_f)=0$ then the deformation space $Def(C',f)$ is smooth, unobstructed, any deformation of the nodes of $C'$ lifts to a deformation of the pair $(C',f)$, and the general deformation of $(C',f)$ is a smooth curve with a map to $S$. Thus, the last claim of the lemma will also follow from the vanishing $H^1(C',\CN_f)=0$. Note that  $H^1(C',\CN_f)=H^1(C',\CN_f/\CN_f^{\rm tor})$ since $H^1(C',\CN_f^{\rm tor})=0$, and let us show that $H^1(C',\CN_f/\CN_f^{\rm tor})=0$.

Since $S$ is smooth, $\Omega_S$ is a locally free sheaf. Consider the exact sequence
$$0\to \CL\to f^*\Omega_S\to \Omega_{C'}\to \CN_f^{\rm tor}\to 0\, .$$
The sheaf $\CL$ is invertible since $f$ is an embedding in a neighborhood of each node of $C'$. Furthermore, $\CN_f^{\rm tor}=\CO_D$ for the divisor $D=\sum_{x\in C'}{\rm ord}_x(df)\cdot x\subset C'$ whose support belongs to the smooth locus of $C'$. Hence the sequence
\begin{equation}\label{eq:exseqL}
0\to \CL\to f^*\Omega_S\to \Omega_{C'}(-D)\to 0
\end{equation}
is exact, and after dualizing it, we obtain an exact sequence
$$0\to \Theta_{C'}(D)\to f^*\Theta_S\to \CL^*\to \CE xt^1_{\CO_{C'}}(\Omega_{C'}(-D), \CO_{C'})\to 0\, .$$

Let $m\subset\CO_{C'}$ be the ideal sheaf of the nodes of $C'$. Then $\CE xt^1_{\CO_{C'}}(\Omega_{C'}, \CO_{C'})\simeq \CE xt^1_{\CO_{C'}}(\Omega_{C'}(-D), \CO_{C'})\simeq \CO_{C'}/m$, and we have an exact sequence
\begin{equation}\label{eq:exseqNfLdual}
0\to \CN_f^{\rm tor}\to \CN_f\to \CL^*\to \CO_{C'}/m\to 0\, .
\end{equation}

Let $\phi\colon \widetilde{C}\to C'$ be the normalization of $C'$. By \eqref{eq:exseqNfLdual}, $\CN_f/\CN_f^{\rm tor}\simeq m\CL^*\simeq\phi_*\CF$ for an invertible sheaf $\CF$ on $\widetilde{C}$, and
$$H^1(C',\CN_f/\CN_f^{\rm tor})=H^1(\widetilde{C}, \CF)=\oplus_{E}H^1(E, \CF|_E),$$ where the sum is taken over all irreducible components $E\subset \widetilde{C}$. Pick $E$, and set $\phi_E:=\phi|_E$. Let $\Delta_E\subset E$ be the preimage of the nodes of $C'$. Then the sequence
$$0\to \CF|_E\to \phi_E^*\CL^*\to \CO_{\Delta_E}\to 0$$
is exact. Hence $c_1(\CF|_E)=c_1(\phi_E^*\CL^*)-|\Delta_E|$.

Since $\CL$ is invertible and $E$ is smooth, the pullback of \eqref{eq:exseqL}
$$0\to \phi_E^*\CL\to \phi_E^*f^*\Omega_S\to \phi_E^*\Omega_{C'}(-D)\to 0$$
is exact. Thus, $c_1(\phi_E^*\CL)=f(\phi_E(E)).K_S-c_1(\Omega_E)-|\Delta_E|+\deg(\phi_E^*D)$, and \begin{equation}\label{eq:c1formula}
c_1(\CF|_E)=-f(\phi_E(E)).K_S+c_1(\Omega_E)-\deg(\phi_E^*D)>c_1(\Omega_E)
\end{equation}
by the assumption of the lemma. Hence $H^1(C',\CN_f)=\oplus_{E}H^1(E, \CF|_E)=0$.

It remains to compute the dimension of $Def(C',f)$. Since $Def(C',f)$ is smooth, $\dim Def(C',f)=\dim Def^1(C',f)$, and by \eqref{eq:seqdef} we obtain
$$\dim Def(C',f)=h^0(C',\CN_f)+\dim \CO_{C'}/m=h^0(C',\CN_f)+p_a(C')-p_g(C')\, .$$
$$
\begin{array}{rl}
  h^0(C',\CN_f)=& h^0(C',\CN_f^{\rm tor})+h^0(\widetilde{C},\CF)\\
  = &\deg(D)+\sum_E \left(c_1(\CF|_E)-p_g(E)+1\right)=-C.K_S+p_g(C')-1
\end{array}
$$
by \eqref{eq:exseqNfLdual}, \eqref{eq:c1formula}, and Riemann-Roch theorem. Hence
$$\dim Def(C',f)=-C.K_S+p_a(C')-1$$
as expected. The proof is now complete.
\end{proof}

\subsubsection{Proof of Theorem~\ref{thm:redSev}}
Proposition~\ref{prop:existnodal} (2), (4), and Proposition~\ref{prop:existmixed}, which we prove below, imply that under the assumptions of Theorem~\ref{thm:redSev}, the Severi varieties contain at least two irreducible components, one of which parameterizes curves having no nodes at all, and another parameterizing curves having at least one node. This completes the proof of Theorem~\ref{thm:redSev}.

\begin{proposition}\label{prop:existnodal}
Let $d\ge 1$, $0\le g_1\le \frac{(d-1)(d-2)}{2}$ and $0\le g_2\le (d-1)^2$ be integers. Then

(1) $V^{irr}(\PP^2,\CO_{\PP^2}(d), g_1)$ contains an irreducible component of maximal dimension $3d+g_1-1$, whose general closed point corresponds to a nodal curve.

(2) Assume that $p>2$. Then the Severi variety $V^{irr}(S_q,\CL^d_q,g_1)$ contains an irreducible component $V$, such that the curve corresponding to a general closed point of $V$ has no nodes among its singularities.

(3) $V^{irr}(\PP^1\times\PP^1,\CO_{\PP^1\times\PP^1}(d,d), g_2)$ contains an irreducible component of maximal dimension $4d+g_2-1$, whose general closed point corresponds to a nodal curve.

(4) The Severi variety $V^{irr}(S'_q,(\CL'_q)^d,g_2)$ contains an irreducible component $V'$, such that the curve corresponding to a general closed point of $V'$ has no nodes among its singularities.
\end{proposition}
\begin{proposition}\label{prop:existmixed} Let $d\ge 2$, $\frac{q-1}{2}\le g_1\le \frac{2dq-2d-q-1}{2}$,  and $q-1\le g_2\le 2dq-q-d-1$ be positive integers.

(1) Assume that $p>2$. Then $V^{irr}(S_q,\CL_q^d, g_1)$ contains an irreducible component of maximal possible dimension $3d+g_1-1$, whose general closed point corresponds to a curve having at least one node.

(2) $V^{irr}(S'_q,(\CL'_q)^d, g_2)$ contains an irreducible component of maximal possible dimension $4d+g_2-1$, whose general closed point corresponds to a curve having at least one node.
\end{proposition}
\begin{proof}[Proof of Proposition~\ref{prop:existnodal}]
(1) Let $C$ be the union of $d$ general lines. Then $C$ has only nodes as its singularities. Mark $\delta=\frac{(d-1)(d-2)}{2}-g_1$ nodes such that their complement in $C$ is connected (e.g. mark only nodes that do not belong to one of the lines), and let $x_1,\dotsc,x_{d+g_1}$ be the unmarked nodes. Let $(C',f)$ be as in Lemma~\ref{lem:defcur}. Then $C'$ is connected, $p_a(C')=g_1$, and $C.K_{\PP^2}=-3d$. Hence the dimension of $Def(C',f)$ is $3d+g_1-1$, and by Lemma~\ref{lem:defcur}, a general deformation of $(C',f)$ consists of an irreducible curve with a map to $\PP^2$. Since $C'\to C$ is birational, the fiber over $C$ of the natural projection $Def(C',f)\to |\CO_{\PP^2}(d)|$ is finite. Thus, $V^{irr}(\PP^2,\CO_{\PP^2}(d), g_1)$ contains an irreducible component $\widetilde{V}$ of dimension $3d+g_1-1$ whose general closed point corresponds to a nodal curve.

(2) By Theorem~\ref{thm:genZar} (1), $$\dim(V^{irr}(S_q,\CL^d_q,g_1))\le -\CL^d_q.K_{S_q}+g_1-1=3d+g_1-1=
\dim(V^{irr}(\PP^2,\CO_{\PP^2}(d),g_1)).$$
Let $\widetilde{V}\subset V^{irr}(\PP^2,\CO_{\PP^2}(d), g_1)$ be an irreducible component as in (1). Recall that $S_q=\PP^2/\mu_q$. Let $\pi\colon \PP^2\to S_q$ be the natural projection, and $D\subset \PP^2$ be a curve corresponding to a closed point of $\widetilde{V}$. Then $\pi(D)$ corresponds to a closed point of $V^{irr}(S_q,\CL^d_q,g_1)$. Since $\pi$ is bijective all the singularities of $\pi(D)\in |\CL^d_q|$ are unibranch, but the images of the nodes of $D$, which are mapped to singularities of type $A_{2q-1}$. The induced map $\widetilde{\pi}\colon \widetilde{V}\to V^{irr}(S_q,\CL^d_q,g_1)$ is injective, and $V:=\widetilde{\pi}(\widetilde{V})$ satisfies the required condition.

The proofs of (3) and (4) are identical to the proofs of (1) and (2). We leave the details to the reader.
\end{proof}

\begin{proof}[Proof of Proposition~\ref{prop:existmixed}]
Let $C$ be the union of a general curve $E\in|\CL_q|$ and $d-1$ curves $C_1,\dotsc, C_{d-1}$ corresponding to general points of $V^{irr}(S_q,\CL_q, 0)$. The curve $\cup_{i=1}^{d-1}C_i$ has no nodes by Theorem~\ref{thm:char>2}. Thus, $C$ has $(d-1)q$ nodes, since $\CL_q.\CL_q=2Area(\Delta_q)=q$ and
$E$ is a general curve in a very ample linear system.

Mark $\frac{2dq-2d-q+1}{2}-g_1$ nodes on $C$ in such a way that for each $i$ at least one of the nodes in $E\cap C_i$ is not marked. Let $x_1,\dotsc,x_k$ be the remaining nodes. Then $k=g_1+\frac{2d-q-1}{2}$. Let $(C',f)$ be as in Lemma~\ref{lem:defcur}. Then $C'$ is connected, $p_a(C')=p_a(E)+k-(d-1)=g_1$, and the irreducible components of $C'$ are: $E$ and the normalizations $\PP^1_i$ of $C_i$. Furthermore, by Theorem~\ref{thm:char>2}, each $\PP^1_i$ contains a unique point at which the differential of the map $\PP^1_i\to \Sigma_q$ vanishes, and the order of vanishing at this point is one, since in local coordinates $f$ is given by $t\mapsto (t^2,t^q)$ (cf. the second paragraph in the proof of Theorem~\ref{thm:char>2}). Thus, the assumptions of Lemma~\ref{lem:defcur} are satisfied, since $K_{\Sigma_q}.C_i=K_{\Sigma_q}.E=-3$.

Hence the dimension of $Def(C',f)$ is $3d+g_1-1$, and a general deformation of $(C',f)$ consists of an irreducible curve with a map to $S_q$ by Lemma~\ref{lem:defcur}. Since $C'\to C$ is birational, the fiber over $C$ of the natural projection $Def(C',f)\to |\CL_q^d|$ is finite. Thus, $V^{irr}(S_q,\CL_q^d, g_1)$ contains an irreducible component of dimension $3d+g_1-1$. Furthermore, since each marked node of $C$ has two preimages in $C'$, the curve corresponding to a general closed point of the component has at least $\frac{2dq-2d-q+1}{2}-g_1\ge 1$ nodes.

The proof of (2) is almost identical. The only difference here is that each rational component contains two points at which the differential $df$ vanishes if $p>2$, and unique such point if $p=2$. In the first case, the order of vanishing of $df$ at each point is one. In the second case it is two, since in local coordinates the map is given by $t\mapsto (\frac{t^2}{t-\lambda},t^q)$, $\lambda\ne 0$. Hence one can use Lemma~\ref{lem:defcur} again to obtain the result. We leave the details to the reader.
\end{proof}

\end{document}